\documentclass{article}
\usepackage{geometry}		     
\usepackage{amssymb,xcolor,xspace}
\usepackage[tworuled,vlined,oldcommands]{algorithm2e}
\usepackage{amsmath,graphicx,comment}

\title{On the genericity of Whitehead minimality\thanks{This work was partially supported by the
ANR through ANR-2010-BLAN-0204, through ANR-10-LABX-58 and through ANR-JCJC-12-JS02-012-01}}

\author{Fr\'ed\'erique Bassino\\
\small{Universit\'e Paris 13, Sorbonne Paris Cit\'e, LIPN, CNRS, (UMR 7030)}\\
\small{F--93430, Villetaneuse, France. \email{bassino@lipn.univ-paris13.fr}}
\and Cyril Nicaud\\
\small{Universit\'e Paris-Est, LIGM, CNRS UMR 8049}\\
\small{F-77454 Marne-la-Vall\'ee, France. \email{nicaud@univ-mlv.fr}}
\and Pascal Weil\\
\small{CNRS, LaBRI, UMR 5800, F-33400 Talence, France. \email{pascal.weil@labri.fr}}\\
\small{Univ. Bordeaux, LaBRI, UMR 5800, F-33400 Talence, France}}

\date{\today}

\def\cqfd{\skip10=\parfillskip\parfillskip=0pt
\enspace\hfill\symbolecqfd\par\parfillskip=\skip10\par\medskip}
\def\symbolecqfd{\rlap{$\sqcap$}$\sqcup$}

\newtheorem{theorem}{Theorem}[section]
\newtheorem{proposition}[theorem]{Proposition}
\newtheorem{lemma}[theorem]{Lemma}
\newtheorem{corollary}[theorem]{Corollary}

\newtheorem{pro-fact}[theorem]{Fact}

\newtheorem{pro-example}[theorem]{Example}
\newenvironment{example}{\begin{pro-example}\rm}{\cqfd\end{pro-example}}

\newtheorem{pro-remark}[theorem]{Remark}
\newenvironment{remark}{\begin{pro-remark}\rm}{\cqfd\end{pro-remark}}

\newenvironment{preuvee}{\rm \trivlist \item[\hskip \labelsep{\bf
Proof.}]}{\cqfd\endtrivlist}

\def\cqfd{\skip10=\parfillskip\parfillskip=0pt
\enspace\hfill\symbolecqfd\par\parfillskip=\skip10\par\medskip}
\def\symbolecqfd{\rlap{$\sqcap$}$\sqcup$}
\def\preuve{\begin{preuvee}}
\def\eop{\end{preuvee}}

\newcounter{commentcounter}

\newenvironment{Pitemize}
  {\begin{list}{}%
      {\setlength{\itemsep}{0pt}%
       \setlength{\topsep}{0pt}%
      }%
  }%
{\end{list}}

\def\binom#1#2{{#1\choose#2}}

\def\para#1{{\medskip\noindent\bf #1.}}

\def\inter[#1]{[\![#1]\!]}
\def\inv{^{-1}}

\def \Iem {S^{k}}
\def \Ilm {S^{\leq k}}
\def \Igm {S^{\geq k}}

\def \P {\mathbb{P}}

\def \O {\mathcal{O}}

\def \D {\mathcal{D}}
\def \calE {\mathcal{E}}
\def \B {\mathcal{B}}

\def \I {\mathcal{I}}

\def \calR {\mathcal{R}}

\let\phi\varphi
\let\epsilon\varepsilon

\def\debut{\operatorname{\tt begin}}

\def\sequence{\operatorname{\tt sequence}}
\def\extr{\operatorname{\tt extr}}
\def\positive{\operatorname{\tt positive}}
\def\negative{\operatorname{\tt negative}}

\def\milieu{\operatorname{\tt mid}}

\def\Wh{\operatorname{\sf Wh}}

\def \oa {\overline{a}}
\def \ob {\overline{b}}

\def\red {\mathcal{R}}
\def\cred {\mathcal{C}}

\def\FB#1 {\textcolor{red}{#1}\xspace}
\def\CN#1 {\textcolor{blue}{#1}\xspace}
\def\PW#1 {\textcolor{purple}{#1}\xspace}

\def\email{\texttt}

\begin{document}

\maketitle

\begin{abstract}
We show that a finitely generated subgroup of a free group, chosen uniformly at random, is strictly Whitehead minimal with overwhelming probability. Whitehead minimality is one of the key elements of the solution of the orbit problem in free groups. The proofs strongly rely on combinatorial tools, notably those of analytic combinatorics. The result we prove actually depends implicitly on the choice of a distribution on finitely generated subgroups, and we establish it for the two distributions which appear in the literature on random subgroups.
\end{abstract}

\section{Introduction}

The problem we consider in this paper is the generic complexity of the \emph{Whitehead minimization problem} for finitely generated subgroups of a free group $F(A)$. Every such subgroup $H$ is a regular subset of $F(A)$ and can be represented uniquely by a finite, edge-labeled graph $\Gamma(H)$ subject to particular constraints, called the Stallings graph of the subgroup; this discrete structure constitutes a natural tool to compute with subgroups, and it also provides a notion of size for $H$: we denote by $|H|$ the number of vertices of $\Gamma(H)$.

A natural equivalence relation on subgroups is provided by the action of the automorphism group of $F(A)$: the subgroups $H$ and $K$ are in the same orbit if $K = \phi(H)$ for some automorphism $\phi$ of $F(A)$ --- that is, $H$ and $K$ are ``the same'' up to a change of basis in the ambient group. The Whitehead minimization problem consists in finding a minimum size element in the orbit of a given finitely generated subgroup $H$. This problem is decidable in polynomial time (Roig, Ventura and Weil \cite{2007RoigVW}, following an early result of Gersten \cite{1984:Gersten}). We refer the readers to~\cite{LyndonSchupp} for the usage of this problem in solving the more general orbit membership problem.

Here we are rather interested in the notion of generic complexity, that is, the complexity of the problem when restricted to a generic set of instances (a set of instances such that an instance of size $n$ sits in it with probability tending to 1 when $n$ tends to infinity; precise definitions are given below). Our main result states that the generic complexity of the Whitehead minimization problem is constant, and more precisely, that the set of Whitehead minimal subgroups is generic (see \cite{2004:MiasnikovMyasnikov} for an early discussion of the generic complexity of this problem, especially in the case of cyclic subgroups).

An implicit element of the discussion of complexity is the notion of size of inputs. In the case of finitely generated subgroups of a free group, we can use either a $k$-tuple ($k$ fixed) of words which are generators of the subgroup $H$ (and the size of the input is the sum of the lengths of these words), or the Stallings graph of $H$ (and the size is $|H|$). These two ways of specifying the subgroup $H$ give closely related worst-case complexities (because of linear inequalities between the two notions of size), but they can give very different generic complexities: it was shown in \cite{BMNVW} that malnormality (an important property of subgroups) is generic if subgroups are specified by a tuple of generators, whereas non-malnormality is generic if subgroups are specified by their Stallings graph. Our results show that Whitehead minimality is generic in both set-ups.


A key ingredient of our proofs is a purely combinatorial characterization of Whitehead minimality in terms of the properties of the graph $\Gamma(H)$ (Proposition~\ref{prop from RVW} below), proved in \cite{2007RoigVW}, which involves counting the edges labeled by certain subsets of the alphabet in and out of each vertex. This is what allows us to turn the algebraic problem into a combinatorial one, which can be tackled with the methods of combinatorics and theoretical computer science.

Interestingly, the reasons why Whitehead minimality is generic when subgroups are specified by their Stallings graph, and why it is generic when subgroups are specified by a $k$-tuple of words, are directly opposite. The Stallings graph of the subgroup generated by a $k$-tuple of words of length at most $n$ generically consists of a small central tree and long loops connecting leaves of the tree, so much of the geometry of the graph is along these long loops, where each vertex is adjacent to only two edges. In contrast, an $n$-vertex Stallings graph generically has many transitions and each vertex is adjacent to a near-full set of edges.

The origins of this work go back to discussions with Armando Martino and Enric Ventura in 2009.

\section{Preliminaries}

Let $r > 1$, let $A$ be a finite $r$-element set and let $F(A)$ be the \emph{free group on $A$}. We can think of $F(A)$ as the set of \emph{reduced} words on the symmetrized alphabet $\tilde A = A \cup \bar{A}$, where $\bar A = \{\bar a \mid a\in A\}$. Recall that a word is reduced if it does not contain occurrences of the words of the form $a\bar a$ or $\bar aa$ ($a\in A$). The operation $x\mapsto \bar x$ is extended to $\tilde A^*$ by letting $\bar{\bar a} = a$ and $\overline{ub} = \bar b\bar u$ for $a\in A$, $b\in \tilde A$ and $u\in \tilde A^*$.

We denote by $[n]$ the set of positive integers less than or equal to $n$, and by $\red_n$ (resp. $\red_{\le n}$) the set of reduced words of length exactly (resp. at most) $n$. A reduced word $u$ is called \emph{cyclically reduced} if $u^2$ is reduced, and we let $\cred_n$ (resp. $\cred_{\le n}$) be the set of cyclically reduced words of length exactly (resp. at most) $n$.

\subsection{Stallings graph of a subgroup}\label{sec: stallings}

It is now classical to represent the finitely generated subgroups of a free group by finite rooted edge-labeled graphs, subject to certain combinatorial constraints. An \emph{$A$-graph} is a finite graph $\Gamma$ whose edges are labeled by elements of $A$. It can be seen also as a transition system on alphabet $\tilde A$, with the convention that every $a$-edge from $p$ to $q$ represents an $a$-transition from $p$ to $q$ and an $\bar a$-transition from $q$ to $p$. Say that $\Gamma$ is \emph{reduced} if
it is connected and if no two edges with the same label start (resp. end) at the same vertex: this is equivalent to stating that the corresponding transition system is deterministic and co-deterministic.  If 1 is a vertex of $\Gamma$, we say that $(\Gamma,1)$ is \emph{rooted} if every vertex, except possibly 1, has valency at least 2.

If $H$ is a finitely generated subgroup of $F(A)$, there exists a unique reduced rooted graph $(\Gamma(H),1)$, called the \emph{Stallings graph} of $H$, such that $H$ is exactly the set of reduced words accepted by $(\Gamma(H),1)$: a reduced word is accepted when it labels a loop starting and ending at $1$. Moreover, this graph can be effectively computed given a tuple of reduced words generating $H$, in time $\O(n\log^*n)$ \cite{1983Stallings,2006Touikan}. We denote by $|H|$ the number of vertices of $\Gamma(H)$, which we interpret as a notion of \emph{size} of $H$. Observe that if $H$ is the cyclic subgroup generated by a cyclically reduced word $w$, then $|H|$ is the length of $w$. This algorithmic construction and the idea of systematically using these graphs to compute with finitely generated subgroups of free groups, go back to Serre's and Stallings' seminal papers (\cite{1977Serre} and \cite{1983Stallings} respectively).

\begin{figure}[htbp]
\begin{minipage}{.24\textwidth}
\centering
\includegraphics[scale=.8]{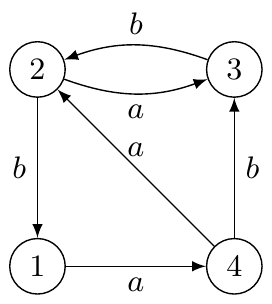}
\end{minipage}
\begin{minipage}{.75\textwidth}
\caption{The Stallings graph of $H = \langle aab, ab\oa b, abbb\rangle$. The reduced word $u=aa\ob\oa b$ is in $H$ as it is accepted by $\Gamma(H)$: it labels a path starting from $1$ and ending at $1$, with edges being used backward when reading a negative letter. Since every vertex has valency at least 2, this graph is cyclically reduced.
\label{fig:stallings}}
\end{minipage}
\end{figure}

We record the following fact, which will be useful in the sequel. Say that an $A$-graph $\Gamma$ is \emph{cyclically reduced} if it is reduced and every vertex has valency at least 2. The $A$-graph in Fig.~\ref{fig:stallings} is cyclically reduced. If $H$ is a finitely generated subgroup of $F(A)$ and $\Gamma(H)$ is not cyclically reduced, then the distinguished vertex 1 has valency 1. Let $\Gamma'$ be the graph obtained from $\Gamma(H)$ by repeatedly erasing every vertex of valency 1 (and the edges adjacent to them): then $\Gamma'$ is cyclically reduced and if $v$ is a vertex of $\Gamma'$, then $(\Gamma',v)$ is the Stallings graph of some conjugate $H^g = g\inv Hg$ of $H$.

\subsection{Whitehead minimality}\label{subsec: whitehead}

Say that a subgroup $H$ is \emph{Whitehead minimal} if it has minimum size in its automorphic orbit, that is if $|H| \le |\phi(H)|$ for every automorphism $\phi$ of $F(A)$. It is \emph{strictly Whitehead minimal} if $|H| < |\phi(H)|$ for every automorphism $\phi$ that is not length preserving (i.e., that is not induced by a permutation of $\tilde A$). Strict Whitehead minimality means that $H$ is the only minimum size representative of its orbit, up to a permutation of the letters (that is, up to a relabeling of the edges of its Stallings graph).

Observe, following the discussion at the end of Section~\ref{sec: stallings}, that if $\Gamma(H)$ is not cyclically reduced, then $H$ is not Whitehead minimal.

A crucial characterization of (strict) Whitehead minimality can be expressed in terms of the so-called \emph{Whitehead automorphisms}. More precisely Whitehead exhibited a finite family $\Wh(A)$ of automorphisms of $F(A)$, with the remarkable property that a subgroup is Whitehead minimal if and only if $|H| \le |\phi(H)|$ for every $\phi\in \Wh(A)$ (this is a result of Whitehead himself for cyclic subgroups, see \cite{LyndonSchupp}, and of Gersten in the general case \cite{1984:Gersten}).

In this paper we will use a combinatorial formulation of this characterization of Whitehead minimality, which was proved in \cite{2007RoigVW}, and which we now explain. We distinguish three kinds of Whitehead automorphisms. Firstly, the length-preserving automorphisms of $F(A)$, which permute the letters of $\tilde A$ and for which we always have $|\phi(H)| = |H|$: they can be disregarded when assessing whether a subgroup is Whitehead minimal. Secondly the inner automorphisms of the form $g\mapsto g^v = v\inv gv$ for some letter $v\in \tilde A$. As discussed above, $\Gamma(H)$ is not cyclically reduced if and only if one of these automorphisms satisfies $|\phi(H)| < |H|$.

The third and last kind of Whitehead automorphisms is in bijection with the set of pairs $(Y,v)$ where $Y$ is a subset of $\tilde A$ and $v$ is a letter in $\tilde A$ such that $v\in Y$, $\bar v\not\in Y$ and
$2 \le |Y| \le 2|A| - 2$. Such a pair $(Y,v)$ is called a \emph{Whitehead descriptor}. The corresponding Whitehead automorphism fixes the letters $v$ and $\bar v$ and maps each letter $a\in \tilde A \setminus \{v,\bar v\}$	to
    $$\phi(a) = v^\lambda a v^\rho \textrm{ where }
    \lambda = \begin{cases}-1 & \textrm{if $\bar a\in Y$,} \cr 0 & \textrm{otherwise;}\end{cases}
    \qquad \rho = \begin{cases}1 & \textrm{if $a\in Y$,} \cr 0 & \textrm{otherwise.}\end{cases}$$
Let $\Gamma$ be a reduced graph, and let $(Y,v)$ be a Whitehead descriptor. Then we let $\positive(\Gamma,Y,v)$ be the set of vertices of $\Gamma$
with at least one incoming edge labeled by a letter in $Y$, at least
one incoming edge labeled by a letter not in $Y$, and no incoming edge
labeled $v$.  Let also $\negative(\Gamma,Y,v)$ be the set of vertices with an incoming edge labeled $v$, and all other incoming
edges labeled by letters in $Y$.

\begin{example}
Consider the Whitehead descriptor $(Y,v)$ with $v=\oa$ and $Y=\{\oa,b\}$. For the graph $\Gamma$ depicted on Fig.~\ref{fig:stallings}, vertex $1$ is in $\negative(\Gamma,Y,v)$ since its incoming edges are labeled by $b$ and $\oa$ (obtained by flipping the edge $1\xrightarrow{a} 4$). Vertex $3$ is in $\positive(\Gamma,Y,v)$ since its incoming edges are labeled by $a$, $b$  and $\ob$, one not in $Y$, one in $Y$ and all different from $v$. One can also verify that vertices $2$ and $4$ are neither in $\positive(\Gamma,Y,v)$ nor in $\negative(\Gamma,Y,v)$.
\end{example}

The following statement is a reformulation of the Whitehead-Gersten characterization of Whitehead minimality mentioned above in terms of these parameters; it is a consequence of \cite[Proposition 2.4]{2007RoigVW}.

\begin{proposition}\label{prop from RVW}
    A finitely generated subgroup $H$ of $F(A)$ is Whitehead
    minimal (resp. strictly Whitehead
    minimal) if and only if it is cyclically reduced and, for every
    Whitehead descriptor $(Y,v)$, we have $|\positive(\Gamma(H),Y,v)| \ge |\negative(\Gamma(H),Y,v)|$ (resp. $|\positive(\Gamma(H),Y,v)| > |\negative(\Gamma(H),Y,v)|$).
\end{proposition}

\preuve
Proposition 2.4 in \cite{2007RoigVW} actually states that, if $(Y,v)$ is a Whitehead descriptor and $\phi$ is the corresponding Whitehead automorphism, then $|\phi(H)| - |H| = |C(H)| - |D(H)|$, where $C(H)$ is the set of vertices of $\Gamma(H)$ with incoming $Y$-labeled and $Y^c$-labeled edges, and $D(H)$ is the set of vertices with an incoming $v$-labeled edge. The intersection $B(H) = C(H) \cap D(H)$ is the set of vertices with an incoming $v$-labeled edge and some incoming $Y^c$-labeled edge. Moreover, $\positive(\Gamma(H),Y,v)$ is the complement of $B(H)$ in $C(H)$ and $\negative(\Gamma(H),Y,v)$ is the complement of $B(H)$ in $D(H)$. The proposition follows immediately.
\eop

\subsection{Distributions over finitely generated subgroups}\label{subsec: proba}

Let $S$ be a countable set, the disjoint union of finite sets $S_n$
($n\ge 0$), and let $B_n = \bigcup_{i\le n}S_i$.  Typically in this
paper, $S$ will be the set of Stallings graphs, of partial injections,
of reduced words or of $k$-tuples of reduced words, and $S_n$ will be
the set of elements of $S$ of size $n$.

A subset $X$ of $S$ is \emph{negligible} if the probability for an element of $B_n$ to be in $X$, tends to 0 when $n$ tends to infinity; that is, if $\lim_n\frac{|X\cap B_n|}{|B_n|} = 0$.

The notion is refined as follows: we say that $X$ is \emph{exponentially} (resp. \emph{super-polynomially}, \emph{polynomially}) \emph{negligible} if $\frac{|X\cap B_n|}{|B_n|}$ is $\O(e^{-cn})$ for some $c > 0$ (resp. $\O(n^{-k})$ for every positive integer $k$, $\O(n^{-k})$ for some positive integer $k$). The set $X$ is exponentially  (resp. super-polynomially, polynomially, simply) \emph{generic} if its complement is exponentially (resp. super-polynomially, polynomially, simply) negligible.
We note the following elementary lemma.

\begin{lemma}\label{lemma: negligible in a subset}
With the above notation, if $C \subseteq S$ satisfies $\liminf_n\frac{|C \cap B_n|}{|B_n|} = p > 0$ and $X$ is exponentially  (resp. super-polynomially, polynomially, simply) negligible in $S$, then so is $X \cap C$ in $C$.
\end{lemma}

\preuve
The verification is immediate if we observe that, for $n$ large enough,
$$\frac{|X \cap C \cap	B_n|}{|C \cap B_n|} \le \frac{|X \cap  B_n|}{|C \cap B_n|} = \frac{|X \cap  B_n|}{|B_n|}\ \frac{|B_n|}{|C \cap B_n|} \le \frac2p\ \frac{|X \cap  B_n|}{|B_n|}.$$
\eop

Genericity and negligibility can also be defined using the radius $n$ spheres $S_n$ instead of the balls $B_n$. The same properties are generic or negligible, exponentially, super-polynomially, polynomially or simply, provided $|B_n|$ grows fast enough, see for instance \cite[Sec. 2.2.2]{BMNVW}.

\para{The graph-based distribution}
The uniform distribution on the set of size $n$ Stallings graphs was
analyzed by Bassino, Nicaud and Weil \cite{BNW}.  Here we summarize
the principles of this distribution and the features which
will be used in this paper.

In a Stallings graph, each letter labels a partial injection on the
vertex set: in fact, such a graph can be viewed as an $A$-tuple
$\vec f = (f_a)_{a\in A}$ of partial injections on an $n$-element set, with a
distinguished vertex, and such that the resulting graph (with an
$a$-labeled edge from $i$ to $j$ if and only if $j = f_a(i)$) is
connected and has no vertex of valency 1, except perhaps the
distinguished vertex.  We may even assume that the $n$-element set in
question is $[n]$, with 1 as the distinguished
vertex, see \cite[Section 1.2]{BNW} for a precise justification.

Let $\I_n$ denote the set of partial injections on $[n]$ and let $\B_n$ be the set of $r$-tuples in $\I_n^r$  which define a Stallings graph (recall that $|A| = r$). Let also $\D_n$ be the subset of $\B_n$, of those $r$-tuples which define a cyclically reduced Stallings graph. Then $\D_n$ (and hence $\B_n$) is generic in $\I_n^r$~\cite[Corollary 2.7]{BNW}

The fundamental observation, used in \cite{BNW} to achieve this result, is the following: the functional graph of a partial injection $f\in \I_n$ (that is: the pair $([n],E)$ where $i\rightarrow j\in E$ whenever $j=f(i)$), is made of cycles and sequences.
%
%
This allows the use of the analytic combinatorics calculus on exponential generating series (EGS) \cite[Sec. II.2]{FlajoletSedgewick}. Recall that, if $I_n$ is the number of partial injections on $[n]$, the corresponding EGS is $I(z) = \sum_{n\geq 0} \frac1{n!}I_n z^n$. From~\cite[Sec. 2.1 and Proposition 2.10]{BNW}, we get
\begin{equation}\label{eq:I(z)}
I(z) = \frac1{1-z}\exp\left(\frac z{1-z}\right) \quad
    \mbox{and} \quad  \frac{I_n}{n!} =\frac{e^{-\frac12}}{2\sqrt\pi}e^{2\sqrt n}n^{-\frac14} (1 + o(1)).
\end{equation}
The formula for $I(z)$ is based on the fact that a partial injection is a set of sequences (whose EGS is $\frac{z}{1-z}$) and of cycles (whose EGS is $\log\big(\frac1{1-z}\big)$). We refer the readers to \cite[Sec. II.2]{FlajoletSedgewick} and \cite{BNW} for further details. We use again this calculus in Section~\ref{sec: statistical properties}.

\para{The word-based distribution}
The distribution more commonly found in the literature (e.g.
\cite{2006KSS,2002Jitsukawa,2003KapovichMSS}), which we term \emph{word-based}, originated in the work of Arzhantseva and Ol'shanski\u\i\ \cite{1996ArjantsevaO}. It is in fact a distribution on the $k$-tuples $\vec h = (h_1, \ldots, h_k)$ of reduced words of length at most $n$, where $k$ is fixed and $n$ is allowed to grow to infinity; one then considers the subgroup $H$ generated by $\vec h$.

This is a reasonable way of defining a distribution on finitely generated subgroups of $F(A)$, and even on rank $k$ subgroups, in spite of the fact that different tuples may generate the same subgroup (see for instance \cite[Sec. 3.1]{BMNVW}).

The literature also considers Gromov's so-called density model, which uses much larger random tuples (of positive density within $\cred_n$). This model is usually considered to study the asymptotic properties of finite group presentations rather than subgroups of $F(A)$ and we will not discuss it here (see for instance \cite{2005:Ollivier}).

We will use the following statistics on the number of reduced and cyclically reduced words, which can be easily verified:
$$|\red_m| = 2r(2r-1)^{m-1}\quad \textrm{and}\quad 2r(2r-1)^{m-2}(2r-2) \le |\cred_m| \le |\red_m|.$$
Summing over all $m \le n$, we find that
$$|\red_{\le n}| = \frac{r}{r-1}\big((2r-1)^n-1\big)\quad \textrm{and}\quad 2r\big((2r-1)^{n-1}-1\big) \le |\cred_{\le n}| \le |\red_{\le n}|.$$
In particular, both $|\red_{\le n}|$ and $|\cred_{\le n}|$ are $\Theta\big((2r-1)^n\big)$ and $\liminf_n\frac{|\cred_{\le n}|}{|\red_{\le n}|} > 0$ (see Lemma~\ref{lemma: negligible in a subset}).

\section{The graph-based distribution}
\label{sec: W minimal graph}

We now study the genericity of strict Whitehead minimality for
the graph-based distribution. The proof of Theorem~\ref{thm: Whitehead graph} below is given in Sections~\ref{sec: statistical properties} and~\ref{sec: from pi to sg}.

\begin{theorem}\label{thm: Whitehead graph}
    Strict Whitehead minimality is super-polynomially generic for the
uniform distribution over the set of cyclically reduced Stallings graphs. 
\end{theorem}
%

%




\subsection{Statistical properties of size $n$ partial injections}\label{sec: statistical properties}

If $f$ is a partial injection on $[n]$, we let
\begin{Pitemize}
\item[$\bullet$] $\sequence(f)$ be the number of sequences in the functional graph of $f$; a sequence has at least one vertex;

\item[$\bullet$]
$\extr(f)=\{i\in[n] \mid f(i)\text{ is undefined or }i\text{ has no preimage by }f\}$; it is
the set of extremities of sequences in the functional graph of $f$.
\end{Pitemize}
%
%
We note that, for every $f\in \I_n$, because of length $1$ sequences,
\begin{equation}\label{eq: seq vs extr}
\sequence(f) \le |\extr(f)| \le 2\sequence(f).
\end{equation}

\begin{proposition}\label{pro:large deviations non-trivial seq}
For the uniform distribution, the probability that the number
of sequences of a size $n$ partial injection is
not in $(\frac12\sqrt{n},2\sqrt{n})$ is super-polynomially small (of the form $\O(e^{-c\sqrt n})$ for some $c>0$).
\end{proposition}

\preuve
If $T(z)$ is a formal power series, we denote by $[z^n]T(z)$ the coefficient of $z^n$ in the series.
For any $k\geq 0$, let $\Iem(z)$, $\Ilm(z)$ and $\Igm(z)$ be the EGSs of the partial
injections having respectively exactly $k$, at most $k$ and at least $k$ sequences. Observe that
an injection with $k$ sequences is a set of $k$ sequences together with a set of cycles; the symbolic method~\cite[Sec. II.2]{FlajoletSedgewick} therefore yields:
\[
\Iem(z) = \frac1{k!}\left(\frac{z}{1-z}\right)^k \frac1{1-z}.
\]
The radius of convergence of this series is 1, and Cauchy's estimate for the coefficient of a power series  \cite[Theorem 10.26]{rudin}  states that for any positive real $\zeta<1$, we have
\[
[z^n]\Iem(z) \leq \frac{\Iem(\zeta)}{\zeta^n}.
\]
Taking $\zeta = 1-\frac1{\sqrt{n}}$ approximatively minimizes the right hand quantity, and after basic
computations we obtain that for $n$ large enough,
\[
[z^n]\Iem(z) \leq \sqrt{n}\;e^{2+\sqrt{n}}\ \frac{n^{\frac{k+1}2}}{k!}.
\]
Since $S^{\leq \frac12\sqrt{n}}(z) = \sum_{k=0}^{\frac12\sqrt{n}} \Iem(z)$ and $S^{\geq 2\sqrt{n}}(z) = \sum_{k=2\sqrt{n}}^{n} \Iem(z)$ we get upper
bounds for coefficients of both series by bounding
$\sum_{k=0}^{\frac12\sqrt{n}} \frac1{k!}n^{\frac{k}2}$ and $\sum_{k=2\sqrt{n}}^{n} \frac1{k!}n^{\frac{k}2}$ from above. The term
$\frac1{k!}n^{\frac{k}2}$ is increasing in the first sum and decreasing in the second one, so we can
bound each term of each series by its maximum value. This yields the following inequalities:
\begin{align*}
\sum_{k=0}^{\frac12\sqrt{n}} \frac{n^{\frac{k}2}}{k!} \leq \sum_{k=0}^{\frac12 \sqrt{n}} \frac{n^{\frac14\sqrt{n}}}{(\frac12\sqrt{n})!},
&\quad [z^n]S^{\leq \frac12\sqrt{n}}(z) \le \frac{n^{\frac32+\frac14\sqrt{n}}}{(\frac12\sqrt{n})!}\;e^{2+\sqrt{n}}
\quad\text{and}\\
\sum_{k=2\sqrt{n}}^{n} \frac{n^{\frac{k}2}}{k!} \leq
\sum_{k=2\sqrt{n}}^{n} \frac{n^{\sqrt{n}}}{(2\sqrt{n})!},
&\quad [z^n]S^{\geq 2\sqrt{n}}(z) \le \frac{n^{2+\sqrt{n}}}{(2\sqrt{n})!}\;e^{2+\sqrt{n}}.
\end{align*}
Using the Stirling bounds \cite[Eq. (9.15), p. 54]{feller} $n!\geq n^ne^{-n}$ and the asymptotics of $I_n$ in
Eq.~\eqref{eq:I(z)}, we obtain upper bounds of the announced form for
\[
\frac{[z^n]S^{\leq \frac12\sqrt{n}}(z)}{[z^n]I(z)}\quad\text{and}\quad\frac{[z^n]S^{\geq 2\sqrt{n}}(z)}{[z^n]I(z)},
\]
respectively the probabilities for a partial injection on $[n]$ to have at most $\frac12\sqrt{n}$ and at least $2\sqrt{n}$ sequences.
\eop

%

We use Proposition~\ref{pro:large deviations non-trivial seq} to bound the number of vertices that are simultaneously extremities for two partial injections.

\begin{proposition}\label{pro:intersection}
    For the uniform distribution over size $n$ pairs of partial injections, the probability
    \[\P \left(|\extr(f) \cap \extr(f')| \ge \frac{\sqrt n}{4(r-1)}\right)\]
    is super-polynomially small (of the form $\O(e^{-c\sqrt n})$ for some $c>0$).
\end{proposition}

\preuve
Let $f$ and $f'$ be partial injection on $[n]$.
By Proposition~\ref{pro:large deviations non-trivial seq} and Eq.~(\ref{eq: seq vs extr}), the probability that one of them has more than
$4\sqrt{n}$ extremities is super-polynomially small --- so we can restrict
the analysis to the cases where both $f$ and $f'$ have at most $4\sqrt{n}$
extremities, up to a super-polynomially small error term.

Let $m=\lfloor 4\sqrt{n}\rfloor$. Let
$E_f$ and $E_{f'}$ be two sets obtained by adding uniformly at
random elements of $[n]$ to $\extr(f)$ and $\extr(f')$ respectively,
until $|E_f|=|E_{f'}|=m$. Note that by symmetry, and since
$f$ and $f'$ are chosen independently, both $E_f$ and $E_{f'}$ are uniform
and independent
size $m$ subsets of $[n]$. Moreover, since $\extr(f)\subseteq E_f$ and
$\extr(f')\subseteq E_{f'}$, we have
\[
\P \left(|\extr(f) \cap \extr(f')| \ge \frac{\sqrt n}{4(r-1)}\right)
\leq \P \left(|E_f \cap E_{f'}| \ge \frac{\sqrt n}{4(r-1)}\right).
\]
It suffices therefore to show that,
super-polynomially generically, the intersection of two $m$-element subsets of $[n]$ has less than
$\frac{\sqrt n}{4(r-1)}$ elements.
Let $X(n,m,k)$ be the number of pairs of $m$-subsets whose intersection has size $k$. Then
\[
X(n,m,k) = \binom{n}{k}\binom{n-k}{m-k}\binom{n-m}{m-k}.
\]
Therefore the probability that the intersection has size $k$ is
\[
\P(|E_{f}\cap E_{f'}| = k) = \frac{X(n,m,k)}{\binom{n}{m}^2} = k!\binom{m}{k}^2 \ \frac{(n-m)!^2}{n!(n-2m+k)!}.
\]
Note that $\frac{(n-m)!^2}{n!(n-2m+k)!} < (n-m)^{-k}$, that $\binom{m}{k} < 2^m$. Let $\alpha=\frac1{4(r-1)}$. Then
$$\P(|E_{f}\cap E_{f'}| \ge \alpha\sqrt n) = \sum_{k=\alpha\sqrt n}^m \P(|E_{f}\cap E_{f'}| = k) < 2^{2m} \sum_{k=\alpha\sqrt n}^m \frac{k!}{(n-m)^k}.$$
Moreover $k\mapsto\frac{k!}{(n-m)^k}$ is decreasing for $k \le m$ (for $n$ large enough), so we have
$$\P(|E_{f}\cap E_{f'}| \ge \alpha\sqrt n) < 2^{2m} m \frac{(\alpha\sqrt n)!}{(n-m)^{\alpha\sqrt n}} < 2^{8\sqrt n} 4 \sqrt n \left(\frac{\alpha\sqrt n}{n-4\sqrt n}\right)^{\alpha\sqrt n}.$$
This concludes the proof since the dominant term is of the form $n^{-\frac\alpha2\sqrt n}$.
\eop

\subsection{From partial injections to Stallings graph}
\label{sec: from pi to sg}

Notice that if $(Y,v)$ is a Whitehead descriptor, the definitions of the functions $\negative(-,Y,v)$ and $\positive(-,Y,v)$ make sense for all $r$-tuple of size $n$ partial injections, even if they do not form a (cyclically reduced) Stallings graph. We will use the following combinatorial bounds to establish Theorem~\ref{thm: Whitehead graph}.

\begin{lemma}\label{lm:bounds neg pos}
Let $(Y,v)$ be a Whitehead descriptor and let $\vec f = (f_a)_{a\in A} \in \I_n^r$. If $v\in \bar A$, we let $f_v = f_{\bar
v}\inv$. Then we have
\begin{align*}
|\negative(\vec{f},Y,v)| & \le	\sum_{a\ne v}|\extr(f_v) \cap \extr(f_a)|, \\
|\positive(\vec{f},Y,v)| & \ge \sequence(f_v) - \sum_{a\ne v}|\extr(f_v) \cap \extr(f_a)|.
\end{align*}
\end{lemma}

\preuve
Recall that a vertex $p$ in $\negative(\vec{f},Y,v)$ has an incoming $v$-edge and all its incoming edges have labels in $Y$. Since $\bar v\not\in Y$, it follows that $p \in \extr(f_v)$. Moreover, if $a\not\in Y$ and $a\ne\bar v$ (there exists such an $a$ since $|Y|\le 2r-2$), $p$ has no incoming $a$-edge, so $p\in\extr(f_a)$. This establishes the first inequality.

Similarly, if $v\in A$ and $p$ is the initial vertex of a sequence of $f_v$ (and hence a $v$-extremity), and if in addition $p$ is not an $a$-extremity for any $a\ne v, \bar v$, then $p\in \positive(\vec{f},Y,v)$. Therefore, if $\debut(f_v)$ denotes the set of initial vertices of sequences of $f_v$, we have
$$\debut(f_v) \setminus \bigcup_{a\ne v, \bar v}\extr(f_v) \cap \extr(f_a) \subseteq  \positive(\vec{f},Y,v),$$
and the announced inequality follows since $|\debut(f_v)| = \sequence(f_v)$.

If $\bar v\in A$ we consider instead the set of final vertices of sequences in $f_{\bar v}$.
\eop

\noindent{\bf Proof of Theorem~\ref{thm: Whitehead graph}}.\enspace
%
%
Let $\D_n$ be the set of $r$-tuples of size $n$ partial injections which define a cyclically reduced Stallings graph, and let $\calE_n$ be the set of $r$-tuples $\vec f$ of size $n$ partial injections which fail to satisfy $|\positive(\vec f,Y,v)|>|\negative(\vec f,Y,v)|$ for some Whitehead descriptor $(Y,v)$. By Proposition~\ref{prop from RVW}, we want to show that $\calE_n \cap \D_n$ is super-polynomially negligible within $\D_n$.

Since $\D_n$ is generic in the full set of $r$-tuples of partial injections, namely $\I_n^r$ (see Section~\ref{subsec: proba}), Lemma~\ref{lemma: negligible in a subset} shows that we only need to show that $\calE_n$ is super-polynomially negligible in $\I_n^r$.

For each Whitehead descriptor $(Y,v)$, let $\calE_n(Y,v)$ denote the set of $r$-tuples $\vec f \in \I_n^r$ such that
$|\positive(\vec f,Y,v)|\leq|\negative(\vec f,Y,v)|$. Then $\calE_n$ is the (finite) union of the $\calE_n(Y,v)$ and it suffices to prove that each $\calE_n(Y,v)$ is super-polynomially negligible in $\I_n^r$.

For a fixed Whitehead descriptor $(Y,v)$, Lemma~\ref{lm:bounds neg pos} shows that
$$\P\Big(\calE_n(Y,v)\Big) \le \P\Big(\sequence(f_v) \leq 2 \sum_{a\ne v}|\extr(f_v) \cap \extr(f_a)|\Big).$$
We observe that if $|\extr(f_v) \cap \extr(f_a)| < \frac{1}{4(r-1)}\sqrt{n}$ for each $a\in A$, $a\ne v,\bar v$ and $\sequence(f_v) > \frac12\sqrt n$, then $2\sum_{a\ne v}|\extr(f_v) \cap \extr(f_a)| < \frac12\sqrt n < \sequence(f_v)$, so that $\vec f \not\in \calE_n(Y,v)$. Therefore, by considering the complements of these properties, we see that $\P(\calE_n(Y,v))$ is at most equal to
\begin{align*}
\P\Big(\sequence(f_v) \le \frac12\sqrt{n}\Big) + \sum_{a\ne v} \P\Big(|\extr(f_v) \cap \extr(f_a)| \ge \frac1{4(r-1)}\sqrt{n}\Big).
\end{align*}
This concludes the proof since each of the summands is super-polynomially small by Propositions~\ref{pro:large deviations non-trivial seq} and~\ref{pro:intersection}.
\cqfd


 Theorem~\ref{thm: Whitehead graph} is stated for the uniform distribution
on \emph{cyclically reduced} Stallings graphs. One may wonder if
a similar result holds for the uniform distribution on Stallings graph. We show the following.

\begin{corollary}
Strict Whitehead minimality is polynomially, but not super-polyn\-omially, generic
for the uniform distribution over Stallings graphs.
\end{corollary}

\preuve
As per the proof of Theorem~\ref{thm: Whitehead graph}, an $r$-tuple $\vec f \in \I_n^r$ satisfies super-polyn\-omially generically the constraint that $|\positive(\vec{f},Y,v)|>|\negative(\vec{f},Y,v)|$ for any Whitehead descriptor $(Y,v)$, -- and hence a Stallings graph $(\Gamma(H),1)$ super-polynomially generically satisfies the constraint $|\positive(\Gamma(H),Y,v)|>|\negative(\Gamma(H),Y,v)|$ for any $(Y,v)$.


For $H$ to be strictly Whitehead minimal, $\Gamma(H)$ must also be cyclically reduced. Equivalently, vertex 1 must be of valency at least 2, that is, it must not be an extremity for one letter and isolated (i.e., the extremity of a length 1 sequence) for all other letters.

The probability that a vertex $p$ is an extremity for the partial injection $f$ is $\frac1n|\extr(f)|$, which is $\Theta(\frac1{\sqrt{n}})$ by Proposition \ref{pro:large deviations non-trivial seq}. The probability that $p$ is isolated is $\frac{I_{n-1}}{I_{n}}$, which is $\Theta(\frac1n)$ by Eq.~(\ref{eq:I(z)}). Therefore, vertex 1 is of valency less than 2 with
probability $\Theta(n^{-(r-1)-\frac12})$, which concludes the proof.
\eop

%

In other words, the uniform distribution on Stallings graphs exhibits the same behavior as that on cyclically reduced graphs with respect to strict Whitehead minimality, but with a weaker error term.

\section{The word-based distribution}
\label{sec: W minimal word}

Let $k\ge 2$ be a fixed integer. We discuss the genericity of strict Whitehead minimality for the subgroups generated by a random $k$-tuple of cyclically reduced words and we show the following.

\begin{theorem}\label{thm:word-based}
For the uniform distribution over $k$-tuples of cyclically reduced words of length at most $n$,
strict Whitehead minimality is exponentially generic.
\end{theorem}

\subsection{Shape of the Stallings graph}

The following elementary statement combines results established in~\cite{1996ArjantsevaO,2002Jitsukawa} and in \cite[Sec. 3.1]{BMNVW}.

\begin{proposition}\label{pro:small heart}
Let $\alpha \in (0,1)$ and $0<\beta<\frac12\alpha$, let $\vec h = (h_1,\ldots, h_k)$ be a tuple of elements of $\red_{\le n}$ and let $H$ be the subgroup generated by $\vec h$. Then, exponentially generically,

- $\min|h_i| > \lceil\alpha n\rceil$ and the prefixes of the $h_i$ and $h_i\inv$ of length $\lfloor\beta n\rfloor$ are pairwise distinct

- the Stallings graph $\Gamma(H)$ consists of a \emph{central tree} of height $\lfloor\beta n\rfloor$ -- whose vertices can be identified with the prefixes and suffixes of length at most $\lfloor\beta n\rfloor$ of the $h_i$ -- and of $k$ \emph{outer loops}, one for each $h_i$, of length $|h_i|-2\lfloor\beta n\rfloor$, connecting the leaves of the central tree.
\end{proposition}

Proposition~\ref{pro:small heart} describes the typical shape of a Stallings graph
under the word-based distribution: as $\beta$ can be taken arbitrarily small and $\alpha$
arbitrarily close to $1$, an overwhelming proportion of the vertices are in the outer loops, and in particular have valency exactly two.

\subsection{Counting the occurrences of short factors}

If $u$ is a word over an alphabet $B$, we denote by $Z_n(u)$ the function that counts the occurrences of $u$ as a factor in a word in $B^n$.

\begin{lemma}\label{lm:uniform random words}
Let $B$ be a finite alphabet with $k\geq 2$ letters and let $u\in B^m$. Then the mean value of $Z_n(u)$ is asymptotically equivalent to $\frac{n}{k^m}$. Moreover, for any $\epsilon>0$ there exists a constant $c>0$ such that
\[
\P\left(\left|Z_n(u)-\frac{n}{k^m}\right|\geq \epsilon n\right) \leq e^{-c n}.
\]
\end{lemma}

\preuve
For $i\in[n+1-m]$, the probability $X_n^{(i)}$ that $u$ is a factor at position $i$ in
a random word of length $n$ is $k^{-m}$, with the convention that the first letter
is at position $1$.
For each $\ell \in [m]$, let $Z^{(\ell)}_n(u) = \sum_j X_n^{(mj+\ell)}$,
for $0\leq j\leq \lfloor\frac{n+1-\ell}{m}\rfloor$.
%
%
Each $Z^{(\ell)}_n(u)$ is the sum of independent random variables since there is no
overlap in the portions of the length $n$ word considered. Therefore $Z^{(\ell)}_n(u)$ follows a binomial law of parameters $k^{-m}$ and $\lfloor\frac{n+1-\ell}{m}\rfloor$: by Hoeffding's inequality \cite{1963:Hoeffding}, it is centered around its mean value which is equivalent to $\frac{n}{mk^m}$, and it satisfies  $\P\left(\left|Z^{(\ell)}_n(u)-\frac{n}{mk^m}\right| > \frac\epsilon mn\right) \le e^{-c_\ell n}$ for some $c_\ell > 0$ and for each $n$ large enough.
The announced result follows from the fact that $Z_n(u) = Z_n^{(0)}(u) + \ldots + Z_n^{(m-1)}(u)$.
\eop

Now if $u$ is a reduced word over the alphabet $\tilde A$, we denote by $\tilde Z_n(u)$ the function that counts the occurrences of $u$ as a factor in a reduced word in $\red_n$.

\begin{lemma}\label{lm:reduced words}
Let $u=u_1u_2$ be a reduced word of length $2$. Then for any $\epsilon>0$ there exists a constant $c>0$ such that, for $n$ large enough,
\[
\P\left(\tilde Z_n(u) > \left(\frac1{(2r-1)^2} + \epsilon\right)(n-1) + 1 \right)  \leq e^{-cn} \]
 and
\[\P\left(\tilde Z_n(u)<\left( \frac{2r-2}{(2r-1)^3}-2\epsilon \right) (n-1)\right)  \leq e^{-cn}
\]
\end{lemma}

\preuve
We first consider the case where $u_1\neq u_2$.
The idea is to use Lemma~\ref{lm:uniform random words} via an encoding of reduced
words. For every $a\in \tilde A$, let $\phi_a$ be a bijective map from
$\tilde{A}\setminus\{{\bar a}\}$ to $[2r-1]$.
Let $\phi$ be the map from the set of reduced words to $\tilde{A}\times [2r-1]^*$ defined
for every reduced word $z = z_1\cdots z_n$ by
\[
\phi(z) = (z_1,\phi_{z_1}(z_2)\phi_{z_2}(z_3)\cdots\phi_{z_{n-1}}(z_n)).
\]
Observe that for every $n > 0$, $\phi$ is a bijection from $\red_n$ to $\tilde{A} \times [2r-1]^{n-1}$, which is computed by an automaton with outputs: the states are the elements of $\tilde A$ and for every $a\in \tilde A$ and $b\ne \bar a$, there is a transition from $a$ to $b$ on input $b$ with output $\phi_a(b)$.
Moreover, the uniform distribution on $\red_n$ is obtained by choosing $z_1$ uniformly in $\tilde A$, $z'$ uniformly in $[2r-1]^{n-1}$, and taking $\phi\inv(z_1,z')$.

We now choose particular functions $\phi_a$: for every $a\neq \bar u_1$, we choose $\phi_a(u_1)=1$. This way every occurrence of $u_1$ (except possibly for the first letter of $z$), is encoded by a $1$ (note that the $1$s provided by $\phi_{\bar u_1}$ do not encode an occurrence of $u_1$). We also require that $\phi_{u_1}(u_2)=2$ and $\phi_a(\bar u_1)=3$ for every $a\neq u_1$: thus every occurrence of $u= u_1u_2$ in $z$ translates to an occurrence of $12$ in $\phi(z)$, and every
occurrence of $\bar u_1$ translates to a $3$ in $\phi(z)$.
See Figure~\ref{fig:encoding} for an example.

\begin{figure}[htbp]
\begin{minipage}{.25\textwidth}
\[
\begin{array}{c|c|c|c|c}
	  & a &\oa& b & \ob \\
\hline
\phi_a	  & 1 & - & 3 & 2   \\
\hline
\phi_{\oa} & - & 3 & 1 & 2   \\
\hline
\phi_b	  & 1 & 3 & 2 & -   \\
\hline
\phi_{\ob} & 1 & 3 & - & 2   \\
\end{array}
\]
\end{minipage}\qquad
\begin{minipage}{.60\textwidth}
\[
{\large
\begin{array}{c|ccccccccccccccc}
z & b & a & \ob & \oa & b & b & b & a & a & \ob & a & b & a & \ob & a\\
\hline
\phi(z)& b & {\bf 1} & {\bf 2}	 & 3   & {\bf 1} & {\bf 2} & 2 & 1 & {\bf 1} & {\bf 2}	 & 1 & 3 & {\bf 1} & {\bf 2}   & 1\\
\end{array}
}
\]
\end{minipage}
\caption{An example of the encoding used in the proof of Lemma~\ref{lm:reduced words}.
The word $z$ above is encoded using the construction associated with the pattern $u=a\ob$:
$a$ is always encoded by a $1$, $\ob$ by a $2$ and the inverse of the first letter, $\oa$, by a $3$. An occurrence of $u$ always corresponds to an occurrence of ${\bf 12}$ in $\phi(z)$, but the opposite is not true: there are false positives, which are always preceded by a $3$. Note also that an occurrence of {\bf 312} does not always correspond to a false positive.
%
\label{fig:encoding}}
\end{figure}

Then for any $t$, we have $\P(\tilde Z_n(u) > t+1) \leq \P(Z_{n-1}(12) > t)$ (the value $t+1$ in the left-hand side of the inequality corresponds to the possibility of an occurrence of $u$ in the leftmost position). For $t = \left(\frac1{(2r-1)^2} + \epsilon\right)(n-1)$, this yields
\begin{align*}
\P\Big(\tilde Z_n(u) > \ (\frac1{(2r-1)^2} + \epsilon)(n-1) + 1\Big)
&\le \P\Big(Z_{n-1}(12) > (\frac1{(2r-1)^2} + \epsilon)(n-1)\Big)\\
& \le \P\Big(|Z_{n-1}(12) - \frac{n-1}{(2r-1)^2}| \ge \epsilon (n-1)\Big).
\end{align*}
The first inequality to be proved then follows from Lemma~\ref{lm:uniform random words} since the pattern $12$ is taken in $[2r-1]^{n-1}$ equipped with the uniform distribution.

Observe that counting occurrences of $12$ overestimates the number of occurrences
of $u$. More specifically, if a false positive occurs, then the said occurrence of $12$ is preceded by a $3$ in $\phi(z)$. Hence, the number of false positives is bounded above by the number of occurrences of $312$ in $\phi(z)$. Therefore $\P(\tilde Z_n(u) < t) \leq \P(Z_{n-1}(12) - Z_{n-1}(312) < t)$. Let then $t = \left(\frac{2r-2}{(2r-1)^3} - 2\epsilon\right)(n-1) = \left(\frac{n-1}{(2r-1)^2} - \epsilon (n-1)\right) - \left(\frac{n-1}{(2r-1)^3} + \epsilon (n-1)\right)$. Then
\begin{align*}
\P\bigg(\tilde Z_n(u) < \ \left(\frac{2r-2}{(2r-1)^3} - 2\epsilon\right)&(n-1) \bigg)\\
&\le \P\left(Z_{n-1}(12) - Z_{n-1}(312) < \left(\frac{2r-2}{(2r-1)^3} - 2\epsilon\right)(n-1)\right)\\
& \le \P\left(|Z_{n-1}(12) - \frac{n-1}{(2r-1)^2}| > \epsilon (n-1)\right) \\
&\qquad\qquad + \P\left(|Z_{n-1}(312) - \frac{n-1}{(2r-1)^3}| > \epsilon (n-1)\right).
\end{align*}
The second inequality to be proved again follows from Lemma~\ref{lm:uniform random words}.

The case $u=u_1u_1$ is handled in the same fashion, except that we have to set $\phi_{u_1}(u_1) = 2$ instead of $1$.
\eop

\begin{remark}
The statement of Lemma~\ref{lm:reduced words}, and even a slighty stronger statement, can also be obtained using the theory of Markov chains: a reduced word can be seen as a
path in a specific Markov chain -- where the set of states is $\tilde A$, and there is a transition from $a$ to $b$ with probability $\frac1{2r-1}$ whenever $a\neq\bar b$. The result in Lemma~\ref{lm:reduced words} then follows from \cite[Thm 1.1]{lezaud98}. We chose instead to give the elementary and self-contained presentation above.
\end{remark}

\subsection{Proof of Theorem~\ref{thm:word-based}}

Let $\alpha\in(0,1)$, $\beta\in(0,\frac{\alpha}2)$ and $\epsilon>0$ be real numbers, to be chosen later.
Let $W_{n,\alpha,\beta}$ be the set of $k$-tuples $\vec h =(h_1,\ldots,h_k)$ of reduced words of length at most $n$, such that $\min|h_i| > \lceil\alpha n\rceil$ and the prefixes of the $h_i$ and $h_i\inv$ of length $\lfloor\beta n\rfloor$ are pairwise distinct.

For each word $h$ of length greater than $2\lfloor\beta n\rfloor$, let $\milieu(h)$ be the factor of $h$ obtained by deleting the length $\lfloor\beta n\rfloor$ prefix and suffix.

Now let $(Y,v)$ be a Whitehead descriptor and let $H$ be the subgroup generated by $\vec h \in W_{n,\alpha,\beta}$. We denote by $Y^c$ the complement of $Y$. The central tree of $\Gamma(H)$ has at most $2k\beta n$ vertices, and the outer loops of $\Gamma(H)$ are labeled by the $\milieu(h_i)$. All the vertices in these loops have valency 2. Any one of these vertices is in $\negative(\Gamma(H),Y,v)$ if and only if it has an incoming $v$-edge and an outgoing $y$-edge for some $y\in Y^c\setminus\{v\}$. Let $N = (Y\bar v \cup v\bar Y)\setminus\{v\bar v\}$. Then the number of negative vertices in the outer loops is equal to the number of occurrences of elements of $N$ as factors in the $\milieu(h_i)$. That is:
$$\negative(\Gamma(H),Y,v) \leq \sum_{i=1}^k\sum_{xy \in  N} \tilde Z_{|\milieu(h_i)|}(xy) + 2k\beta n.$$

By Proposition~\ref{pro:small heart}, $W_{n,\alpha,\beta}$ is exponentially generic. Moreover, the map $h \mapsto \milieu(h)$ turns the uniform distribution on words in $\red_\ell$ ($\ell > \alpha n$) into the uniform distribution on $\red_{\ell - 2\lfloor\beta n\rfloor}$: indeed, if $u \in \red_{\ell - 2\lfloor\beta n\rfloor}$, then $\P(\milieu(h) = u) = (2r-1)^{-2\lfloor\beta n\rfloor}$, which does not depend on $u$. It follows that the same map also turns the uniform distribution on the set of reduced words of length greater than $\alpha n$ and less than or equal to $n$, into the uniform distribution on its image. Therefore, exponentially generically, we have
\begin{align*}
\negative(\Gamma(H),Y,v) &\leq 2k\beta n + k|N|\Big(\big(\frac{1}{(2r-1)^2}+\epsilon\big)(1-2\beta)n+1\Big)\\
&\le 2k\beta n + 2k(|Y|-1)\Big(\big(1-2\beta\big)\big(\frac{1}{(2r-1)^2}+\epsilon\big)n + 1\Big).
\end{align*}
Similarly, a loop vertex is in $\positive(\Gamma(H),Y,v)$ if it has an incoming $x$-edge with $x \in Y\setminus\{v\}$ and an outgoing $y$-edge with $\bar y\in Y^c$: if $P = (Y\setminus\{v\})\overline{Y^c} \cup Y^c(\bar Y\setminus\{\bar v\})$, then the number of positive vertices in the outer loops is equal to the number of occurrences of elements of $P$ as factors in the $\milieu(h_i)$. That is, exponentially generically,
\begin{align*}
\positive(\Gamma(H),Y,v) &\geq \sum_{i=1}^k\sum_{xy \in  P} Z_{|\milieu(h_i)|}(xy) \\
& \geq k|P|\left(\frac{2r-2}{(2r-1)^3}-2\epsilon\right)\left((\alpha-2\beta)n-1\right)\\
& \geq 2k(|Y|-1)(2r-|Y|)\left(\frac{2r-2}{(2r-1)^3}-2\epsilon\right)\left((\alpha-2\beta)n-1\right).
\end{align*}
In order to conclude, we only need to show that we can choose $\alpha$, $\beta$ and $\epsilon$ such that

\begin{align*}
(2r-|Y|)\Big(&\frac{2r-2}{(2r-1)^3}-2\epsilon\Big)\left((\alpha-2\beta)n-1\right)\\
&> \left(1-2\beta\right)\Big(\frac{1}{(2r-1)^2}+\epsilon\Big) n +1  + \frac{\beta n}{|Y|-1}.
\end{align*}
for all $n$ large enough. The first term is $\Theta(\gamma n)$ with $\gamma = (2r-|Y|)(\frac{2r-2}{(2r-1)^3}-2\epsilon)(\alpha-2\beta)$ and the second term is $\Theta(\delta n)$ with $\delta = (1-2\beta)(\frac{1}{(2r-1)^2}+\epsilon)+\frac{\beta}{|Y|-1}$, so we need to select $\alpha$, $\beta$ and $\epsilon$ such that $\gamma > \delta$. This is possible by continuity, since the limits of these two quantities when $(\alpha,\beta,\epsilon)$ tends to $(1,0,0)$ are respectively $(2r-|Y|)\frac{2r-2}{(2r-1)^3}$ and $\frac1{(2r-1)^2}$, and we have $2r-|Y| \ge 2$ and $\frac{2r-2}{2r-1}\ge \frac23$, so that $(2r-|Y|)\frac{2r-2}{(2r-1)^3} \ge \frac43\ \frac1{(2r-1)^2}$.

This establishes that if $H$ is generated by a $k$-tuple of reduced words, then exponentially generically $\positive(\Gamma(H),Y,v) > \negative(\Gamma(H),Y,v)$ for each Whitehead descriptor. The same exponential genericity holds for $k$-tuples of cyclically reduced words in view of Lemma~\ref{lemma: negligible in a subset} and the discussion at the end of Section~\ref{subsec: proba}. Together with Proposition~\ref{prop from RVW}, this concludes the proof since a subgroup generated by a tuple of cyclically reduced words has a cyclically reduced Stallings graph.
\cqfd\smallskip

To complete the picture, we observe that given a random $k$-tuple of reduced words, instead of cyclically reduced words, there is a non-negligible probability that the graph is not cyclically reduced.

\begin{proposition}\label{lm:non-cyclically-reduced}
For the uniform distribution over $k$-tuples of reduced words of length at most $n$ the Stallings graph
is  not generically cyclically reduced.
\end{proposition}
%

\preuve
Let $ \vec h =(h_1,\ldots,h_k)$ be a random $k$-tuple of reduced
words of length at most $n$ and let $\Gamma(H)$ be the Stalling graph of
the subgroup $H$ generated by $ \vec h$.

We show that with probability tending to
$(\frac1{2r})^{2k-1}$, $\Gamma(H)$ is not cyclically reduced and, more precisely, there exists a letter $a\in\tilde{A}$ such that
every $h_i$ starts with $a$ and ends with $\oa$.

For every pair of letters $a$ and $b$ in $\tilde{A}$,
let $\calR_{a,b}$ be the set of reduced words that start
with  $a$ and end by $b$. Let $R_{a,b}(z)$ be the (ordinary)
generating series associated with $\calR_{a,b}$ defined by
\[
R_{a,b}(z) = \sum_{u\in\calR_{a,b}} z^{|u|}.
\]
Assume that $b\notin\{a,\bar a\}$. Since a word
of $\calR_{a,b}$ is either $ab$ or a word in some
$\calR_{a,c}$ ($c\neq \bar b$) followed
by $b$, we have
\[
R_{a,b}(z) = z^{2} + \sum_{c\neq\bar b}R_{a,c}(z)z,
\]
and similarly
\[
R_{a,a}(z) = z^{2} + \sum_{c\neq\bar a}R_{a,c}(z)z
\qquad\textrm{and}\qquad
R_{a,\bar a}(z) = \sum_{c\neq a}R_{a,c}(z)z.
\]
Now observe that if $b, c \in \tilde{A}\setminus\{a,\bar a\}$, then $R_{a,b}(z) = R_{a,c}(z)$ by symmetry. Hence, fixing a letter $b \in \tilde{A}\setminus\{a,\bar a\}$, the equations above rewrite as
\[
\begin{cases}
R_{a,b}(z) &= z^{2} + (2r-3)R_{a,b}(z)z + R_{a,a}(z)z + R_{a,\bar a}(z)z\\
R_{a,a}(z) &= z^{2} + (2r-2)R_{a,b}(z)z+ R_{a,a}(z)z \\
R_{a,\bar a}(z) &= (2r-2)R_{a,b}(z)z + R_{a,\bar a}(z)z.
\end{cases}
\]
Solving this system yields (thank you \texttt{maple}!)
\begin{align*}
R_{a,\bar a}(z) & = \frac{2z^{3}(r-1)}{(1-z^{2})(1-(2r-1)z)}\\
& = \frac{2r-2}{2r-1} - \frac{1}{2(1-z)} - \frac{r-1}{2r(1+z)}
+ \frac{1}{2r(2r-1)(1-(2r-1)z)}.
\end{align*}
It follows that the number of  words of length $n$ in $\calR_{a,\bar a}$ is
asymptotically equivalent to $\frac{1}{2r}(2r-1)^{n-1}$, and the probability
that a reduced word of length $n$ begins with $a$ and ends with $\bar a$
is asymptotically  equivalent to $\frac{1}{(2r)^2}$. This result also holds
for words of length at most $n$, as they are generically of length greater than
$\frac12 n$.

Thus the probability that the $k$-words of $\vec h$
all begin with the same letter $a$ and end with $\bar a$  is asymptotivally equivalent to
$\frac{1}{(2r)^{2k}}$, and the probability that they all begin with the same letter and end with its opposite
is equivalent to $\frac{1}{(2r)^{2k-1}}$,
which concludes the proof.
\eop

\section{Application to random generation}
%


Proposition~\ref{prop from RVW} and the fact that there are finitely many Whitehead descriptors immediately yield algorithms \texttt{MinimalityTest} (resp. \texttt{StrictMinimalityTest}) to test whether $H$ is (strictly) Whitehead minimal: it suffices to verify whether $\Gamma(H)$ is cyclically reduced (in time at most linear) and to compute,  for each Whitehead descriptor $(Y,v)$, $|\positive(\Gamma(H),Y,v)|$ and $|\negative(\Gamma(H),Y,v)|$. The time required is linear in $|H|$ for each $(Y,v)$, but the number of Whitehead descriptors is exponential in $A$: the resulting algorithm is linear in $|H|$ but not in $|A|$.

In this section, our purpose is different: we want to design efficient random generators -- in the graph-based or the word-based distribution -- for the Stallings graphs of subgroups that are (strictly) Whitehead minimal.

Our algorithms will be rejection algorithms. In general, suppose that $S$ is a countable set, $S$ is the disjoint union of the $S_n$, and $C \subseteq S$ is such that $\liminf_n \frac{|C\cap B_n|}{|B_n]} = p > 0$ (see Section~\ref{subsec: proba} and Lemma~\ref{lemma: negligible in a subset}). If \texttt{RandomS} is a random generator for elements of $S$ and $\texttt{TestC}$ is an algorithm to test whether an element of $S$ is in $C$, then the algorithm in Figure~\ref{fig:randomalgorithm} is a random generator for elements of $C$.

\begin{figure}[htbp]%
    \linesnumbered 
    \centering
    \begin{algorithm}[H]%
      \dontprintsemicolon
      \nocaptionofalgo
      keep $\leftarrow$ False\;
      \Repeat{\text{\emph{keep == True}}} {
	$x$ = \texttt{RandomS}($n$)\;
	keep $\leftarrow$ \texttt{TestC}($x$)\;
	  }
      {\bf return} $x$\;
      \caption{{\tt RandomC}($n$)}
    \end{algorithm}
    \caption{An algorithm to randomly generate an element of $C$ of size $n$}
\label{fig:randomalgorithm}
\end{figure}

In such an algorithm, the loop (lines 3--4) is performed in average $\frac1p$ times. in particular, if both \texttt{RandomS} and \texttt{TestC} take linear time in average, then so does \texttt{RandomC}.

A random generator {\tt RandomStallingsGraph} working in linear average time, is available for the graph-based and the word-based distributions.
\begin{itemize}
\item For the graph-based distribution, such an algorithm is given in~\cite{BNW}.
\item For the word-based distribution, one first generates a $k$-tuple of reduced words (in linear time); next one applies Touikan's algorithm~\cite{2006Touikan} to compute the associated
Stallings graph; it was noted in~\cite[Theorem 4.1]{2012:BassinoNicaudWeil} that the average time complexity of this algorithm is linear.
\end{itemize}

Following the model of the algorithm in Figure~\ref{fig:randomalgorithm}, a rejection algorithm to randomly generate Whitehead minimal subgroups is shown in Figure~\ref{fig: strict algorithm}.

\begin{figure}[htbp]%
    \linesnumbered 
    \centering
    \begin{algorithm}[H]%
      \dontprintsemicolon
      \nocaptionofalgo
      keep $\leftarrow$ False\;
      \Repeat{\text{\emph{keep == True}}} {
	$\Gamma$ = \texttt{RandomStallingsGraph}($n$,$A$)\;
	keep $\leftarrow$ \texttt{MinimalityTest}($\Gamma$)\;
	  }
      {\bf return} $\Gamma$\;
      \caption{{\tt RandomWhiteheadMinimalGraph}($n$,$A$)}
    \end{algorithm}
    \caption{An algorithm to randomly generate Whitehead minimal subgroups}
\label{fig: strict algorithm}
\end{figure}

Similarly, an algorithm {\tt RandomStrictlyWhiteheadMinimalGraph} to randomly generate strictly Whitehead minimal subgroups, is obtained by replacing the call to \texttt{MinimalityTest} by a call to \texttt{StrictMinimalityTest}. In view of the discussion at the beginning of this section, this yields the following statement.

\begin{proposition}\label{pro: strict minimality complexity}
For the graph-based and the word-based distributions,
the average time complexity of the algorithms {\tt RandomWhiteheadMinimalGraph} and {\tt RandomStrictlyWhiteheadMinimalGraph} is linear.
\end{proposition}

\bibliographystyle{plain}

%
%

\end{document}